\documentclass[12pt]{amsart}
\usepackage{amsmath,amssymb}

\newtheorem{thm}{Theorem}%[section]
\newtheorem{lem}[thm]{Lemma}%[section]
\newtheorem{prop}[thm]{Proposition}%[section]
\newtheorem{cor}[thm]{Corollary}%[section]
\theoremstyle{definition}

\theoremstyle{remark}
\theoremstyle{plain}

\newcommand{\Z}{{\mathbf Z}}
\newcommand{\Q}{{\mathbf Q}}

\newcommand{\C}{{\mathbf C}}

\newcommand{\tr}{\operatorname{tr}}

\newcommand{\CN}{{\mathcal N}}

\newcommand{\TT}{\mathbf T^2}

\newcommand{\OPN}{\operatorname{Op}_N}

\newcommand{\HN}{\mathcal H_N}
\newcommand{\TN}{T_N}  %{\hat T_N}
\newcommand{\UN}{U_N}
\newcommand{\Fh}{\mathcal F_h} % h-Fourier transform

\newcommand{\Mat}{\operatorname{Mat}} % matrices

\newcommand{\A}{A}
\newcommand{\D}{D}
\newcommand{\T}{\TN}

\newcommand{\OF}{{\mathfrak O}}
\newcommand{\OO}{\OF_K}
\newcommand{\OK}{\OO}

\newcommand{\trace}{\tr}
\newcommand{\mnorm}{ {\mathcal N}} % norm map

\newcommand{\Torus}{{\mathcal C}_A}

\newcommand{\unit}{\epsilon}

\newcommand{\ord}{\operatorname{ord}}

\newcommand{\ds}{D_A} % discriminant = 4(tr(A)^2-4)
\newcommand{\Li}{\operatorname{Li}}

\newcommand{\lcm}{\operatorname{lcm}}
\newcommand{\LL}{\mathcal L}

\newcommand{\V}{n}
\newcommand{\weil}{U}

\newcommand{\simgeq}{\gtrsim}% 
\newcommand{\simleq}{\lesssim}
\newcommand{\N}{{\mathbf N}}

\numberwithin{equation}{section}

\begin{document}

\title[On quantum ergodicity...]
{On quantum ergodicity for linear maps of the torus}
\author{P\"ar Kurlberg and Ze\'ev Rudnick}
\address{Raymond and Beverly Sackler School of Mathematical Sciences,
Tel Aviv University, Tel Aviv 69978, Israel.  
Current address: Department of Mathematics, 
University of Georgia, Athens, GA 30602, U.S.A.
({\tt kurlberg@math.uga.edu})}
\address{Raymond and Beverly Sackler School of Mathematical Sciences,
Tel Aviv University, Tel Aviv 69978, Israel 
({\tt rudnick@math.tau.ac.il})}

\date{September 27, 1999}
\thanks{Supported in part by a grant from the Israel Science
  Foundation.}

\begin{abstract}

We prove a strong version of quantum ergodicity for linear hyperbolic
maps of the torus (``cat maps''). We show that there is a density one
sequence of integers so that as $N$ tends to infinity along this
sequence, all eigenfunctions of the quantum
propagator at inverse Planck constant $N$ are uniformly distributed. 

A key step in the argument is to show that for a hyperbolic matrix in
the modular group, 
there is a density one sequence  of integers $N$ for which its order
(or period) modulo $N$ is somewhat larger than $\sqrt{N}$. 

\end{abstract}

\maketitle

\section{Introduction}

\subsection{Quantum ergodicity}

An important model for understanding the quantization of classically
chaotic systems are {\em quantum maps}, and in particular the
quantizations of linear automorphisms of the torus $\TT$ 
(``cat maps'').   
Iterating such a  map, we get a discrete dynamical system, 
well-known to be chaotic if the map is hyperbolic. 
A quantization of these ``cat maps'' was proposed by
Hannay and Berry \cite{BH}, see also \cite{Knabe, DE, DEGI}.  
In brief, this procedure restricts Planck's constant to be
an inverse integer: $h=1/N$, and the Hilbert space of states $\HN$ is
$N$-dimensional, in keeping with the intuition that each  state
occupies a Planck cell of volume $h = 1/N$ and the constraint that
the total phase-space $\TT$ has volume one. 
Classical observables (i.e. functions $f\in C^\infty(\TT)$) give rise
to operators $\OPN(f)$ on $\HN$.  
Given a linear automorphism $A$ of the torus, its quantization is a unitary
operator $\UN(A)$ on $\HN$, called the quantum propagator, or
``quantized cat map''.  The eigenfunctions of $\UN(A)$ play the
r\^{o}le of energy eigenstates.

In this paper we will use the quantized cat map to illuminate one of
the few rigorous results available on the semi-classical limit of
eigenstates of classically chaotic systems, namely {\em Quantum
Ergodicity} \cite{Sch, CdV, Zelditch87}. 
%The semi-classical interpretation of this notion was formulated as
%follows: 
To formulate this notion, recall that 
if the classical dynamics are {\em ergodic}, 
then almost all trajectories of a  particle cover the energy
shell uniformly. The intuition 
afforded by the ``Correspondence Principle'' leads one to look for an 
analogous statement about the semi-classical limit of 
expectation values of observables in an energy eigenstate. 
As formulated  by Schnirelman \cite{Sch}, the corresponding assertion
is that when the classical dynamics is {\em ergodic}, 
for almost all eigenstates
the expectation values of observables converge to the phase-space average.
For quantum maps, the form that this takes is the following 
(\cite{Bouzouina96, Zelditch C^*, Zelditch97}): 
Fix an observable $f\in C^\infty(\TT)$. 
%\begin{thm}[Quantum ergodicity for cat maps]\label{qe for cats} 
%Let $A$ be a cat map. 
Then  for  any orthonormal basis $\psi_j$ of $\HN$ consisting of 
eigenfunctions of $\UN(A)$, there is a subset $J(N)\subset
\{1,2,\dots,N\} $, with 
$\frac{\#J(N)}{N}\to 1$, so that for $j\in J(N)$ we have: 
\begin{equation}\label{equdistribution}
\langle \OPN(f)\psi_j,\psi_j \rangle \to \int_{\TT} f,\qquad 
\text{as }N\to \infty
\end{equation}
%\end{thm} 
This is a consequence, using positivity and a standard diagonalization 
argument, %which we omit
of the following estimate for  the variance due to Zelditch
\cite{Zelditch C^*}: Given $f\in C^\infty(\TT)$, 
%\begin{thm} [\cite{Zelditch C^*}]\label{schnir-zel:thm}
%Let $A$ be a  cat map. 
for any orthonormal basis $\psi_j$, $j=1,\dots,N$ of of $\HN$ consisting of 
eigenfunctions of $\UN(A)$,  we have 
%that for all observables $f\inC^\infty(\TT)$, 
\begin{equation}\label{schnir-zel:thm}
 \frac 1N\sum_{j=1}^N \left| \langle \OPN(f)\psi_j,\psi_j
\rangle - \int_{\TT} f \right|^2 \to 0
\end{equation}
%\end{thm} 
Note that the result \eqref{schnir-zel:thm} does not guarantee that
{\em all} eigenfunctions in $\HN$ are equidistributed, even for one
single value of $N$.

\subsection{Beyond quantum ergodicity}
In recent work \cite{KR}, we have found that there is a
commutative group of unitary operators  on the state-space which
commute with the quantized map and therefore act on its eigenspaces. 
We called these ``Hecke operators'', in analogy with the setting of
the modular surface. 
We showed that the joint eigenfunctions of these and of
$\UN(A)$ (which we called ``Hecke eigenfunctions'')  are all 
equidistributed, that is  
\eqref{equdistribution} holds for any choice of Hecke eigenfunctions
in $\HN$. 

Not all eigenfunctions of $\UN(A)$ are Hecke eigenfunctions. 
In fact, the Hecke eigenspaces have small dimension (at most $O(\log
\log N)$), while the 
eigenspaces of $\UN(A)$ may have large dimension. In fact, the {\em mean} 
degeneracy is $N/\ord(A,N)$ where $\ord(A,N)$ the {\em order} (or
period) of $A$ modulo $N$, that is the least integer $k\geq 1$  for
which $A^k=I\mod N$. It can be shown (see section~\ref{smallorder}) 
that the mean degeneracy can be
as large as  $N/\log N$ for arbitrarily large $N$. 
However, it is reasonable to expect that {\em all} 
eigenfunctions become equidistributed - that is we have 
quantum {\em unique} ergodicity.

In this paper, we show ergodicity of {\em all}
eigenfunctions of  $\UN(A)$ for almost all integers $N$: 
\begin{thm} \label{main thm}
Let $A$ be a fixed cat map. 
There is a set of integers $\CN^*$ of density one 
so that all eigenfunctions of $\UN(A)$ are equi-distributed, as
$N\to\infty$, $N\in \CN^*$. 
\end{thm}
Previously, the only result giving an
infinite set of $N$ for which all eigenfunctions of $\UN(A)$ become
equi-distributed is by Degli-Esposti, Graffi and Isola \cite{DEGI}, 
which conditional on  GRH give an infinite set of {\em primes}.

%This result is of a very different nature
%than the one given by \eqref{schnir-zel:thm}, and in
%particular is deterministic. 

\subsection{Outline of the argument} 
Our main tool in relating this result to more traditional themes of 
Number Theory is the following estimate %(Corollary~\ref{}) 
for the fourth power moment of the expectation values,  
giving a condition in terms  of the order of $A$ modulo $N$:  
%that is the least integer $k\geq 1$  for which $A^k=I\mod N$: 
\begin{thm}\label{fourth moment}
There is a sequence of integers of
density one so that for all observables $f\in C^\infty(\TT)$ and 
any orthonormal basis $\{\psi_j\}_{j=1}^N$ of $\HN$ consisting of
eigenfunctions of $\UN(A)$  we have: 
%we have the following estimate for the fourth power moment of the
%expectation values: 
$$
\sum_{j=1}^N |\langle \OPN(f)\psi_j,\psi_j \rangle -\int_{\TT} f |^4 \ll
\frac{N(\log N)^{14}}{\ord(A,N)^{2}} \,.
$$
\end{thm}
%By taking an orthonormal basis with $\psi_1=\psi$ and omitting all but
%the first term we recover Theorem~\ref{quanthm}. 

Thus for any subsequence of integers $N$ such that 
\begin{equation}\label{condition}
\frac{\ord(A,N)}{N^{1/2}(\log N)^{7}}\to\infty
\end{equation}
(and satisfying an additional ``genericity'' assumption explained in
section~\ref{tracetrick:s}) we find that for all eigenfunctions of $\UN(A)$, 
$\langle \OPN(f)\psi,\psi \rangle \to \int_{\TT} f$ as $N\to \infty$.

%It is natural to conjecture the following quantitative form of Quantum Unique
%Ergodicity: For all $N$ and any normalized eigenfunction $\psi$ we have ?
%$$
%|\langle \OPN(f)\psi,\psi \rangle -\int_{\TT} f | \ll
%\frac{1}{\ord(A,N)^{1/2}} 
%$$

Theorem~\ref{fourth moment} reduces the problem of quantum ergodicity to
that of finding sequences of integers satisfying \eqref{condition}, a
problem closely related to the classical Gauss-Artin problem of
showing that any integer, other than $\pm 1$ or a perfect square, is a
primitive root modulo infinitely many primes.  
We show (Theorem~\ref{t:large-order-for-most-N}) that there is some
$\delta>0$ for which there is  
a  set of integers of density $1$ so that 
$$
\ord(A,N)\gg N^{1/2}\exp((\log N)^\delta) \, .
$$
This,  combined with Theorem~\ref{fourth moment} 
gives Theorem~\ref{main thm}. 

To prove Theorem~\ref{t:large-order-for-most-N}, 
we first show in Section~\ref{sec:L(N)} that on a set of density one,
$\ord(A,N)$ is not much smaller than the product of the orders of $A$
modulo prime divisors of $N$. 
Next, we deal with
prime values of $N$.  
In Section~\ref{Goldfeld section} we  show  
(Theorem~\ref{order:thm})  
%\begin{prop}\label{prop ord p} 
that given $1/2<\eta<3/5$, there is a set of primes of positive density
$c(\eta)>0$ so that  
$\ord(A,p) \gg p^{\eta}$. 
We note that this is far short of the truth; by invoking the
Generalized Riemann Hypothesis, one can show that 
for a set of primes of density one, we have $\ord(A,p)\gg p/\log p$ 
(c.f. \cite{erdos-murty}). 
In Section~\ref{ae section} 
we prove Theorem~\ref{t:large-order-for-most-N} by using
Theorem~\ref{order:thm} together with the elementary observation that
for almost all primes $p$, $\ord(A,p)\geq p^{1/2}/\log p$.

As is apparent from this discussion, our result hinge on the condition
\eqref{condition} being satisfied; we can say nothing for $N$ for
which this condition fails, of which there are infinitely many
examples. We consider it a fundamental problem to get results 
when $\ord(A,N)$ is smaller than $N^{1/2}$.

\subsection{Notation}  
We will use the standard convention of analytic number theory: 
Thus $e(z)$ stands for $e^{2\pi iz}$, 
$f(x)\ll g(x)$ as $x\to \infty$ means that there is some $C>0$ so that 
for $x$ sufficiently large, $f(x)<Cg(x)$. Similarly, 
$f(x)\simleq g(x)$ as $x\to \infty$ means  $\limsup f(x)/g(x) \leq
1$. 
We will write $p^t||n$ if $p^t$ divides $n$ but $p^{t+1}$ does not. 
We will denote by $\omega(N)$ the number of prime divisors of $N$.

\newpage
\section{Quantum mechanics on the torus}

\subsection{The Hilbert space of states} 
We review the basics of quantum mechanics on the torus $\TT$, viewed as a
phase space  \cite{BH, Knabe, DE, DEGI}, beginning with a description of 
the Hilbert space of states of such a system: 
We take  state vectors to be
distributions on the line which are periodic in both momentum and
position representations: $\psi(q+1)=\psi(q)$,
$[\Fh\psi](p+1)=[\Fh\psi](p)$, where 
$[\Fh\psi](p) = h^{-1/2}\int\psi(q)\,e(-pq/h)\,dq$.  
The space of such distributions is finite
dimensional, of dimension precisely $N=1/h$, and consists of
periodic point-masses at the coordinates $q=Q/N$, $Q\in \Z$. 
We may then identify $\HN$ with the 
$N$-dimensional vector space $L^2(\Z/N\Z)$, with
the inner product $\langle\,\cdot\,,\,\cdot\,\rangle$ defined by
\begin{equation}
\langle \phi,\psi \rangle 
 = \frac1N \sum_{Q\bmod N} \phi(Q) \, \overline\psi(Q) ,
\end{equation}
%This inner product induces a norm $\norm{\cdot}$ on the space of operators
%on $\HN$, that is on the space of $N\times N$ matrices.

\subsection{Observables} 
Next we construct  quantum observables: 
A central role is played by the translation operators
$$
[t_1 \psi](Q) = \psi(Q+1)
$$
and
$$
[t_2 \psi](Q) = e_N(Q)\, \psi(Q),
$$
which may be viewed as the analogues of differentiation and
multiplication (respectively) operators. 
In fact in terms of the usual translation operators on the line 
$\hat q \psi(q)=q\psi(q)$ and 
$\hat p\psi(q)=\frac{h}{2\pi i} \frac{d}{dq}\psi(q)$, they are given by    
$t_1=e(\hat p)$, $t_2=e(\hat q)$. 
In this context, Heisenberg's commutation relations read 
\begin{equation} \label{commrel}
t_1^a t_2^b = t_2^b t_1^a e_N(ab)  \qquad \forall a,b\in\Z.
\end{equation}

More generally, mixed translation operators are defined for
$n=(n_1,n_2)\in\Z^2$ by  
$$
T_N(n) = e_N(\frac{n_1 n_2}{2}) t_2^{n_2} t_1^{n_1} .
$$
These are unitary operators on $\HN$, whose action on a wave-function
$\psi\in L^2(\Z/N\Z)$ is  given by:  
\begin{equation}\label{action of T(n)}
\TN(n)\psi(Q) = e^{\frac {i\pi n_1 n_2}N} e(\frac{n_2Q}N)\psi(Q+n_1)
\end{equation}
The adjoint/inverse of $\TN(n)$ is given by 
\begin{equation}\label{T(n)^*}
\TN(n)^* = \TN(-n) \,.
\end{equation}
As follows from the commutation relation \eqref{commrel}, we have 
\begin{equation}\label{commute-rel}
T_N(m)\,T_N(n)=e_N(\frac{\omega(m,n)}{2}) \, T_N(m+n)
\end{equation}
where $\omega(m,n)$ is the symplectic form
$$
\omega(m,n)=m_1 n_2 - m_2 n_1 .
$$

For any smooth function $f\in\C^\infty(\TT)$,  define a {\em quantum
observable} $\OPN(f)$, called the {\em Weyl quantization} of $f$ 
\cite{Folland} 
$$
\OPN(f) = \sum_{n\in\Z^2} \widehat f(n) T_N(n)
$$
where $\widehat f(n)$ are the Fourier coefficients of $f$.

Given a state $\psi\in \HN$, the {\em expectation value} of the observable
$f$ in the state $\psi$ is defined to be 
$\langle \OPN(f)\psi,\psi \rangle$.

%We recall some basic properties of the quantized observables 
%$\OPN(f)$: By \eqref{T(n)^*}, the adjoint is given by 
%\begin{equation}\label{adjoint}
%\OPN(f)^* = \OPN(\bar f)
%\end{equation}
%We also need the following fact on the composition of operators: 
%For $f,g\in C^\infty(\TT)$, 
%\begin{equation}\label{composition}
%\OPN(f)\OPN(g) = \OPN(fg)+O_{f,g}(\frac 1N)
%\end{equation}

%Another fact is an analogue of Weyl's law: 
%\begin{equation}\label{weyl}
%\frac 1N \tr \OPN(f) = \int_{\TT} f +O_f(\frac 1N)
%\end{equation}
%This is a consequence of $\tr \TN(0) = N$ and otherwise
%$$
%|\tr\TN(n)|=\begin{cases} N, & n=0\mod N\\0,&\text{otherwise}
%             \end{cases} 
%$$
%In particular for {\em any} orthonormal basis of $\HN$ we get 
%$$
%\frac 1N\sum_{j=1}^N \langle \OPN(f)\psi_j,\psi_j \rangle = \int_{\TT} f
%+O_f(\frac 1N)
%$$
%Note that this is a purely {\em kinematic} fact, independant of any
%dynamics. 

\subsection{Cat maps} 
To introduce dynamics, we consider a  linear automorphism of
the torus $A\in SL(2,\Z)$. 
The iteration of $A$ gives a (discrete) dynamical system, well-known
to be chaotic %\cite{?} 
if $A$ is hyperbolic, that is $|\tr A|>2$
(such a map is called a ``cat map'' in the physics literature).  

%An important quantity for us will be the period, or order, of $A$
%modulo $N$, denoted by $\ord_N(A)$. 
%This is the least integer $k\geq 0$ so that $A^k=I\mod N$.   
%A basic fact is that it always divides the function $\Phi(N) =
%N\prod_{p\mid N}(1-\chi(p))$.  
%In particular,  for all $\epsilon>0$ we have
%$$
%\ord_N(A)\ll N^{2+\epsilon}
%$$ 
%(elaborate on this paragraph...)

If we further assume $A$ is
``quantizable'', that is 
$A=\begin{pmatrix} a&b\\c&d \end{pmatrix}$ with $ab\equiv cd\equiv 0
\mod 2$, then on can assign to $A$ a unitary operator $\UN(A)$ on
$\HN$, the {\em quantum propagator},  whose iterates give the
evolution of the quantum system, and characterized by the property 
(an analogue of ``Egorov's theorem''): 
\begin{equation}\label{Egorov}
\UN(A)^*\OPN(f)\UN(A) = \OPN(f\circ A)
\end{equation}
%Therefore by \eqref{compare Weyl &AW} we have 
%$$
%\UN(A)^*\sOPN(f)\UN(A) = \sOPN(f\circ A)+O_f(\frac 1N)
%$$
This can be thought of as saying that the evolution of the
quantum observable  $\OPN(f)$ follows the evolution 
$f\mapsto f\circ A$ of the classical observable $f$.  
%\footnote{An advantage of Weyl quantization 
%is that for the cat maps, Egorov's theorem is {\em exact}.} 
That \eqref{Egorov} holds {\em exactly} is a special feature of the
linearity of the map $A$; for general maps, \eqref{Egorov} is only
expected to hold asymptotically as $N\to \infty$ (c.f. \cite{MR}). 

The stationary states of the quantum system are given by the
eigenfunctions $\psi$ of $\UN(A)$. It is our goal to study the limiting 
expectation values $\langle \OPN(f)\psi,\psi \rangle$ of observables
in (normalized) eigenstates and show that outside a zero density
set of $N$'s, they all converge to the classical average $\int_{\TT} f$
of the observable as $N\to\infty$. 

\newpage 
\section{The order of a matrix modulo $N$}

\subsection{}
Let $A\in SL(2,\Z)$ be a hyperbolic matrix, that is $|\tr(A)|>2$. 
The {\em order} (or  period)  $\ord(A,N)$ of the map $A$ modulo $N$
is the least integer $k\geq 0$ so that $A^k=I\mod N$. 
We begin to study the order of $A$ modulo an arbitrary integer $N$, 
starting with some well-known generalities. 

\subsubsection{} 
Firstly if $M$ and $N$ are co-prime then 
$$
\ord(A,MN)=\lcm (\ord(A,M),\ord(A,N) )
$$
and so if $N$ has  a prime factorization $N=\prod p_i^{k_i}$ then 
$$
\ord(A,N) = \lcm\{\ord(A,p_i^{k_i})\}
$$

%\subsection{Integral matrices and quadratic fields}
%\label{s:cft}

\subsubsection{}\label{corresponedence}  
The eigenvalues $\unit,\unit^{-1}$ of $A$ generate a field extension 
$K= \Q(\unit)$, which is a real quadratic field since
$\trace(A)^2>4$. We label them so that $|\unit|>1$. 
Let 
$$
D_A = 4(\tr(A)^2-4) 
$$ 
so that  $K=\Q(\sqrt{D_A})$.
We denote by $\OK$ the ring of integers of $K$. 
The eigenvalues $\unit,\unit^{-1}$ of $A$ will be units in $\OK$. 
Adjoining $\unit$ to $\Z$ gives an {\em order} $\OF=\Z[\unit]\subseteq
\OK$ in $K$.  
Then there is an $\OF$-ideal $I\subset\OF$  so
that the action of $\unit$  by multiplication on $I$ is equivalent to
the action of $A$ on $\Z^2$, in the sense that there is  a basis of $I$ with
respect to which the matrix of $\unit$ is precisely $A$ (see
\cite{Taussky} or \cite{KR}). 
The action of $\OF$ by multiplication on $I$ gives us 
an embedding  
$$
\iota:\OF\hookrightarrow \Mat_2(\Z)
$$ 
so that $\gamma = x+y\unit \in \OF$ corresponds to $xI+yA$. Moreover, the
determinant of $xI+yA$ equals $\mnorm(\gamma)=\gamma \bar\gamma$, 
where $\mnorm : K \rightarrow \Q$ is the Galois norm. 
%$\det(\iota(\gamma)) =\mnorm(\gamma)$. 
In particular, if $\gamma \in \OF$ 
has norm one then $\gamma$ corresponds to an element in $SL_2(\Z)$

Given an integer $N\geq 1$, the
embedding $\iota: \OF \hookrightarrow \Mat_2(\Z)$ induces a map
$\iota_N:\OF/N\OF\to \Mat_2(\Z/N\Z)$  and 
the norm $\mnorm:K\to \Q$ gives a
well-defined map 
$$
\mnorm : \OF/N\OF \rightarrow \Z/N\Z . 
$$ 
Denote by $\Torus(N)$ the group of norm one elements in
$\OF/N\OF$: %  (this is an algebraic group defined over $\Z$):  
$$ 
\Torus(N) = \ker \left[ \mnorm:(\OF/N\OF)^*\to(\Z/N\Z)^* \right]\;.
$$
This is a subgroup of $SL(2,\Z/N\Z)$, containing the residues class of
$A$ modulo $N$. 

The cardinality of $\Torus(N)$ can be computed via the Chinese
Remainder Theorem from the cardinality at prime power arguments. 
To do so, define 
$$
\chi(p) =\begin{cases} +1, & p \text{ splits in }K\\
 -1,& p  \text{ inert in }K. \end{cases}
$$
By quadratic reciprocity, $\chi$ is a Dirichlet
character modulo $D_A$  (not necessarily primitive). 
It can then be shown (see e.g. \cite{KR}, Appendix B) 
that if $p$ does not divide $D_A$, then 
\begin{equation}\label{torus(p^k)}
\#\Torus(p^k) = p^{k-1}(p-\chi(p))
\end{equation}
while for primes dividing $D_A$, there is some $c_A>0$ so that 
\begin{equation}\label{ramified torus}
\#\Torus(p^k) \leq c_A p^k \;.
\end{equation}
As a consequence, we find that if $p$ does not divide $D_A$, then the
order of $A$ modulo $p$ divides $p-\chi(p)$, and more generally, for
any prime power $p^k$, if $p$ does not divide $D_A$ then $\ord(A,p^k)$
divides $p^{k-1}(p-\chi(p))$. 

\subsubsection{An upper bound for $\ord(A,N)$} 
Another consequence of \eqref{torus(p^k)}, \eqref{ramified torus} 
is that for any integer $N=\prod p^{k_p}$, 
\begin{equation*}
\#\Torus(N)=\prod_p\#\Torus(p^{k_p}) 
\ll_A N\prod_{p\mid N}(1+\frac 1p)
 \ll_A N \log\log N \;.
\end{equation*}
% If $p\mid D_A$ then $\ord(A,p^k)$ is
%bounded by $c_A p^k$. For any integer $N$, $\ord(A,N)\ll
%bounded by $c_A p^k$. 
Thus, for any integer $N$, we have as an upper
bound for the order 
\begin{equation}\label{no. Torus(N)}
\ord(A,N)\ll N \log \log N\;.
\end{equation}
%N^{1+\epsilon}$ for all $\epsilon>0$.
% changed N^\epsilon to \log\log N. Quote Ramanujan???

\subsection{Making $\ord(A,N)$ small}\label{smallorder}
As for lower bounds on the order, it is easily seen that 
$\ord(A,N) \gg \log N$ for all $N$. 
In fact, this bound is sharp,  as we claim 
\begin{prop}
There is an infinite sequence of
integers $\{N_k\}_{k=1}^\infty$ for which $\ord(A,N_k) \ll \log N_k$. 
\end{prop}
\begin{proof}
To explain the idea, recall first how to find integers $n$ for which
$2$ has small order modulo $n$: The trick is to take $n_k=2^k-1$,
since then $2^k=1\mod n_k$, and so 
$\ord(2,n_k) \leq k \sim \log n_k/\log 2$. 
To modifiy this idea to our context, assume for simplicity 
that the matrix $A$ is ``principal'', that is
the action of $A$ on $\Z^2$ is equivalent to the action of the
unit $\unit$ on the maximal order $\OK$ (in general we need an ideal in
the order $\OF=\Z[\unit]$, see section~\ref{corresponedence}).  
Then $A^k=I \mod N$ is {\em equivalent}  to $\unit^k=1 \mod N\OK$ (in
general, only the implication $\Rightarrow$ is valid). 

Factor $|\det(A^k-I)|$ as a product of prime powers: 
$$
|\det(A^k-I)| = \prod_S p^{\sigma_p} \prod_I p^{\iota_p} \prod_R
p^{\rho_p}
$$
where $\prod_S$ means the product over primes $p=\frak p \bar\frak p$
which split in $K=\Q(\unit)$, $\prod_I$ the product over inert primes
and $\prod_R$ the product over the ramified primes $p=\frak p^2$. 

On the other hand, we have 
$$
\det(A^k-I) = \mnorm(\unit^k-1) =-\unit^{-k}(\unit^k-1)^2 \;.
$$
Write the ideal factorization of $\frak a_k:=(\unit^k-1)\OK  $ as 
$$
\frak a_k = \prod_S \frak p^{s'_p}\bar\frak p^{s''_p} \prod_I p^{i_p}
\prod_R \frak p^{r_p} \;.
$$
Since $\frak a_k^2 = \det(A^k-I)\OK$, we get on comparing the prime
exponents that 
$$
2s'_p = 2s''_p = \sigma_p,\qquad \iota_p=2i_p,\qquad \rho_p=r_p \;.
$$
Since $\sigma_p$ is even, we can set 
$$N_k := \prod_S p^{\sigma_p/2}\prod_I p^{i_p} \prod_R p^{[r_p/2]}\;.$$ 
Then 
$$
N_k\leq |\det(A^k-I)| \leq N_k^2 \delta
$$
where $\delta = \prod_Rp$ is  the product of all ramified primes of
$K$. 

We have $\frak a_k\subseteq N_k\OK$ and so 
$\unit^k=1 \mod N_k\OK$, equivalently $A^k=I \mod N_k$. 
Thus we find 
$$
\ord(A,N_k)\leq k\sim\frac{\log|\det(A^k-I)|}{\log\unit} \leq
\frac{\log N_k^2\delta}{\log\unit}  = \frac 2{\log\unit} \log N_k +O(1)
$$
and so $\ord(A,N_k)\ll \log N_k$ as required. 
\end{proof}

%Let $k$ be odd and relatively prime to $D_K$, the discriminant of $K=
%\Q(\unit)$, and put $N_k'=\det(A^k-I)$. 
%
%If $p|D_K$ then the order of $\unit$ modulo $p^k$ divides $2p^{k-1}$.
% should I explain this???
% \unit is close to +-1 since \unit and \unit^{-1} are ``p-close'',
% so +- \unit is in (1+P) etc etc.
%Since $(k,2D_K)=1$ we thus have 
%$|\unit^{\pm k}-1|_p = |\unit^{\pm  1}-1|_p$, 
%where $|\cdot|_p$ denotes the unique absolute value on $K
%\otimes \Q_p$ that extends the absolute value on $\Q_p$. Since
%$\det(A^k-I)= (\unit^k-1)(\unit^{-k}-1)$ we thus find that if $p|D_K$
%and $p^t||N_k'$ then in fact $p^t|\det(A-I)$.

%Thus, with $N_k$ being the largest divisor of $N_k'$ that is coprime
%to $D_K$ we find that $N_k'/N_k$ is bounded. Moreover, since only
%unramified primes divide $N_k$, we find that $\ord(A,N_k) \leq k$,
%and since $\log N_k \sim k \log\unit$ we are done.

\newpage
\section{Large order of $\A$ implies equidistribution}\label{tracetrick:s}

\subsection{}  
In this section we give a relation between the order of the map $A$
modulo $N$ and the distribution of the eigenfunctions of the
quantization $\UN(A)$. 
We start by relating the fourth power-moment of the 
expectation values 
$\langle \T(\V) \psi_i,  \psi_i \rangle$, 
for $\psi_i$ ranging over an orthonormal basis of 
$\weil_N(\A)$-eigenfunctions, to the number of solutions of a certain
equation modulo $N$. 

\begin{prop}
\label{p:basic-bound}
  Let $\{\psi_i\}_{i=1}^N$  be an orthonormal basis of eigenfunctions of
  $\weil_N(\A)$. Then
%If $\V \not \equiv 0 \mod N$ then
\begin{equation}
\label{e:basic-bound}
\sum_{i=1}^N 
|\langle \T(\V) \psi_i,  \psi_i \rangle |^4  
\leq \frac{N}{\ord(A,N)^4}  \nu(N,n)
\end{equation}
where $\nu(N,n)$ is the number of solutions of the congruence 
$$
\V ( \A^i-\A^j + \A^k - \A^l) \equiv 0 \mod N, \qquad   
 1 \leq i,j,k,l \leq \ord(A,N) 
$$
\end{prop}
\begin{proof}
Let
$$
\D(\V) = 
\frac{1}{\ord(A,N)}
\sum_{i=1}^{\ord(A,N)}
\T( \V \A^i ),
$$
and let $t_{ij}$ be the matrix coefficients of $\T(\V)$ expressed
in terms of the basis $\{ \psi_i \}_{i=1}^N$. From \eqref{Egorov}  we
have that  
$$\T( \V \A^i )=
\weil_N^*(\A) \T( \V ) \weil_N^*(\A)
$$
and by assumption $\weil_N(\A)
\psi_i = \lambda_i \psi_i$ for $\lambda_i$ a root of unity. Thus, if
$\{D_{ij}\}_{i,j=1}^N$ are the matrix coefficients of $\D$ in terms of
the basis $\{ \psi_i \}_{i=1}^N$, then
%as in \cite{KR} we get that
\begin{equation}
\label{e:d-chi}
D_{ij} =
\begin{cases}
  t_{ij} & \text{if $\lambda_i =  \lambda_j$,} \\
  0 & \text{otherwise.}
\end{cases}
\end{equation}
If we denote by $\{v_i\}_{i=1}^N$  the column vectors of $\D$, then the
$(k,k)$-entry of $(\D^*\D)^2$ is
$$ ( (\D^*\D)^2 )_{kk} = 
\sum_i \langle v_i,v_k\rangle \langle v_k,v_i \rangle = 
\sum_i |\langle v_i,v_k\rangle |^2,
$$ 
and since $|\langle v_k,v_k \rangle |=\sum_i |\D_{ki}|^2$ we get
$$
\sum_{ \lambda_i =  \lambda_j }  
|t_{ij}|^4
\leq
\trace( (\D^*\D)^2).
$$
Substituting the definition of $\D$ and using \eqref {T(n)^*} and
\eqref{commute-rel},  
%$$
%\T(n)^*=\T(-n), \quad \T(n)\T(m) = \gamma_{m,n}\T(m+n) \quad 
%\text{(where $|\gamma_{m,n}|=1$)}
%$$
we see that $(\D^{*}\D)^2 $ is given by $\ord(A,N)^{-4}$ times a
sum, ranging over $1 \leq i,j,k,l \leq \ord(A,N)$, of terms
\begin{multline*}
\T( \V \A^i ) 
\T( -\V \A^j )
\T( \V \A^k )
\T( -\V \A^l )
=
\gamma_{i,j,k,l}
\T( \V (\A^i  - \A^j + \A^k + \A^l ) )
\end{multline*}
where $\gamma_{i,j,k,l}$ has absolute value one. Now take the trace;
as follows from \eqref{action of T(n)} (see Lemma~4 in \cite{KR}),
the absolute value of the trace of 
$\T(n)$ equals $N$ if $n \equiv (0,0) \mod N$, zero otherwise. The
result now follows by taking absolute values and summing over all
$i,j,k,l$. (For more details, see section 6.2 in \cite{KR}.)
\end{proof}

%In section ?? we will show that for almost all $N$ the following holds:
%\begin{multline}
%\# \{ 1 \leq i,j,k,l \leq \ord(A,N) \ : \ \V ( \A^i-\A^j + \A^k -
%\A^l) \equiv 0 \mod N  \} = \\
% O(\ord(A,N)^{2+\epsilon})
%\end{multline}
%and we obtain the following
%\begin{cor}
%If $\psi$ is a norm one $\weil(\A_N)$-eigenvector, then for almost all
%$N$, 
%\begin{equation}
%|< \psi, \T(\V) \psi >|  \ll 
%\left( \frac{N}{\ord(A,N)^{2-\epsilon}} \right)^{1/4}
%\end{equation}
%\end{cor}

\subsection{A counting problem}

In order to make use of Proposition~\ref{p:basic-bound} we must bound the
number of solutions to
$$
\V ( \A^i-\A^j + \A^k - \A^l) \equiv 0 \mod N, \quad 1 \leq i,j,k,l
\leq \ord(A,N).  
$$
We will show that there are essentially only {\em trivial} solutions
of this equation, i.e. 
\begin{equation*}
(\A^i,\A^k)=(\A^j,\A^l), \quad (\A^i,\A^k)=(\A^l,\A^k), \quad
\text{or } (\A^i,\A^j)=(-\A^k,-\A^l),
\end{equation*}
where the third possibility only happens if there exists $t$ such that
$\A^t=-I$. In terms of the exponents $i,j,k,l$ this means that
\begin{equation}
\label{e:trivial}
(i,k)=(j,l), \quad (i,k)=(l,k), \quad \text{or } (i,j)=(t-k,t-l),
\end{equation}
where equality is to be interpreted as equality modulo the order of $\A$.

\subsubsection{The prime case}
Here we assume $N=p$ is prime.
\begin{lem}
\label{l:count-prime}
Assume that $\V \A$ and $\V$ are linearly independent modulo $p$, and
that the eigenvalues of $\A$ are distinct modulo $p$. Then there are
at most $3 \ord(A,p)^2$ solutions of  
\begin{equation}
\label{e:count-prime}
\V ( \A^i-\A^j + \A^k - \A^l) \equiv 0 \mod p, \ 
1 \leq i,j,k,l \leq \ord(A,p)
\end{equation}
\end{lem}
\begin{proof}
Let $K$ be the real quadratic field containing the eigenvalues of $A$,
and let $K_p$ be the residue class field at the prime $p$, i.e., $K_p=\OK/P$
where $P$ is a prime of $K$ lying above $p$.  $K_p$ has cardinality
$p$ if $p$ splits in $K$, or $p^2$ if $p$ is inert. 
We may diagonalize the reduction of $\A$ modulo $p$ over the field
$K_p$. In the eigenvector basis we have 
$\A' = \begin{pmatrix} \unit & 0 \\ 0 & \unit^{-1} \end{pmatrix}$ and
$\V'=(n_1',n_2')$, where the assumption of linear independence modulo $p$
implies that both $n_1',n_2' \neq 0$ (in $K_p$.) Thus (\ref{e:count-prime}) is
equivalent to the following two equations over $K_p$:
\begin{equation}
\begin{split}
\unit^i-\unit^j + \unit^k - \unit^l &= 0 \\
\unit^{-i}-\unit^{-j} + \unit^{-k} -\unit^{-l} &= 0
\end{split}
\end{equation}
which in turn (see lemma 15 in \cite{KR}) is equivalent to
\begin{equation}
\begin{split}
  \unit^l = \unit^i-\unit^j + \unit^k  \\
  (\unit^k - \unit^i) (\unit^k - \unit^j) (\unit^i+\unit^j) = 0
\end{split}
\end{equation}
Hence $l$ is determined by the triple $(i,j,k)$. Dividing by $\unit^k$
and letting $i'=i-k$ and $j'=j-k$ we rewrite the second equation as
\begin{equation}
\label{e:product-is-zero}
(1 - \unit^{i'}) (1 - \unit^{j'}) (\unit^{i'}+\unit^{j'}) = 0, \quad 1
\leq i',j' \leq \ord(A,p)
\end{equation}
If the first (or second) factor equals zero then $\ord(A,p) \mid i'$
(or $j'$) since the order of $\unit$ in $K_p^{\times}$ equals
$\ord(A,p)$. If the third factor is zero then $\ord(A,p) \mid i'-j'-t$
where $\unit^t = -1$. In each case this leaves $\ord(A,p)$
possibilities for the pair $(i',j')$, and since $k$ is unconstrained
%whereas $l$ is determined by $(i,j,l)$, 
the total number of solutions 
is at most $3 \ord(A,p)^2$.
\end{proof}

{\bf Remark:} The condition of linear independence mod $p$ in
lemma~\ref{l:count-prime} is satisfied for all but finitely many
primes. In fact, if we let 
$$
M = \begin{pmatrix} n_1 & n_2 \\ m_1 & m_2 \end{pmatrix}
$$
where $\V = (n_1,n_2)$ and $\V \A = (m_1, m_2)$, then the
condition of linear dependence is equivalent to $p \mid \det M$.
Now $\det M$ is a {\em nonzero} integer, because $A$ 
has no rational eigenvectors.  
We also note that if the independence condition is not
satisfied then trivially there are at most $\ord(A,p)^4$ solutions to
(\ref{e:count-prime}).

%\subsubsection{The composite case}

\begin{lem} 
\label{l:count-composite-one}
Let $N$ be square free and coprime to $D_A = 4(\tr(A)^2-4)$. 
Assume further that $\V \A$ and $\V$ are linearly
independent modulo $p$ for all $p\mid N$. Then there are at most
$3^{\omega(N)} \ord(A,N)^2$ solutions of
\begin{equation}
\label{e:count-composite-one}
\V ( \A^i-\A^j + \A^k - \A^l) \equiv 0 \mod N, \quad 
1 \leq i,j,k,l \leq \ord(A,N). 
\end{equation}
\end{lem}

\begin{proof}
Let $(i,j,k,l)$ be a solution to (\ref{e:count-composite-one}). If
$p\mid N$
then (\ref{e:count-composite-one}) holds with $N$ replaced by $p$.
Arguing as in lemma~\ref{l:count-prime} one of the three factors in
(\ref{e:product-is-zero}) must be zero, and the vanishing factor
determines which one of the three equations in (\ref{e:trivial}) that
$(i,j,k,l)$ must satisfy modulo $\ord(A,p)$.  For example, if the 
first factor in (\ref{e:product-is-zero}) is zero, then $(i,j) \equiv
(k,l) \mod \ord(A,p)$.

Now, the group generated by $\A$ modulo $N$ is cyclic and isomorphic
to $\oplus_{q \in Q} \Z/q^{a_q}\Z$ where the $q$'s are {\em distinct}
primes. We will denote the $\Z/q^{a_q}\Z$ component of $i$ by $i_q$
and similarly for $j,k,l$. Since $\ord(A,N)$ is equal to the least
common multiple of $\{\ord(A,p)\}_{p\mid N}$ there exists for each $q \in
Q$ at least one prime $p\mid N$ such that $q^{a_q} \parallel \ord(A,p)$.

Claim: if $(i,j,k,l)$ is a solution to (\ref{e:count-composite-one})
then $(i_q,j_q,k_q,l_q)$ satisfies one of the equations in
(\ref{e:trivial}). The reason is as follows: there is a prime $p\mid N$
such that $q^{a_q} \parallel \ord(A,p)$, 
%and from the proof of lemma~\ref{l:count-prime} we know that 
thus one of the equations in (\ref{e:trivial}) is satisfied modulo
$\ord(A,p)$. Since $q^{a_q} \parallel \ord(A,p)$ this implies that
$(i_q,j_q,k_q,l_q)$ satisfies one of the equations in
(\ref{e:trivial}). (Note in particular that this leaves $q^{2a_q}$
possibilities for $(i_q,j_q,k_q,l_q)$ if we specify one of the
equations in (\ref{e:trivial}) to be satisfied).  Now, to each $p\mid N$
there are 3 different types of trivial solutions, and since
$(i_q,j_q,k_q,l_q)$ must satisfy one of the possibilities in
(\ref{e:trivial}) for all $q \in Q$, we obtain that there are at most
$$
3^{\omega(N)} \prod_{q \in Q} q^{2a_q} = 
3^{\omega(N)} \ord(A,N)^2
$$
solutions to (\ref{e:count-composite-one}). 
\end{proof}

In our applications the hypothesis of linear independence might not
hold for all $p\mid N$. However, we have the following
\begin{lem} 
\label{l:count-composite}
Let $N$ be square free. Then there are at most
$$
O_{\A}\left( |n|_2^{8+\epsilon} 3^{\omega(N)} \ord(A,N)^2 \right)
$$
solutions to  
\begin{equation}
\label{e:count-composite-two}
\V ( \A^i-\A^j + \A^k - \A^l) \equiv 0 \mod N, \quad 
1 \leq i,j,k,l \leq \ord(A,N). 
\end{equation}
\end{lem}
\begin{proof} 
By the remark after lemma~\ref{l:count-prime}, linear dependence
modulo $p$ holds if and only if $ p | \det M$, where $|\det M| \ll_{\A}
|n|_2^2$. Let 
$$
N' = \frac N{\gcd(D_A\det M,N)}  \;.
$$  
Then the hypothesis in lemma~\ref{l:count-composite-one} is 
satisfied for $N'$, leaving  
$3^{\omega(N')}\ord(A,N')^2$ possible values for $(i,j,k,l)$ 
modulo $\ord(A,N')$. Now, an element in $\Z/\ord(A,N') \Z$
has exactly $\frac{\ord(A,N)}{\ord(A,N')}$  preimages in
$\Z \cap [1,\ord(A,N)]$. Hence there are at most 
$$
3^{\omega(N')}\ord(A,N')^2 \left( \frac{\ord(A,N)}{\ord(A,N')} \right)^4
$$
solutions to (\ref{e:count-composite-two}). Since
$$
|\det(M)| \ll_{\A} |n|_2^2
$$ 
we get that 
$$
\frac{N}{N'} = \gcd(D_A \det M,N)\leq D_A \det M\ll_{A} |n|_2^2. 
$$
Finally noting that since $N$ is square-free, 
$$
\ord(A,N) =\lcm \left(\ord(A,N'),\ord(A,N/N') \right) \leq \ord(A,N')\cdot
\ord(A,N/N') 
$$
we find (by \eqref{no. Torus(N)}) that 
$$
\frac{\ord(A,N)}{\ord(A,N')} \leq 
\ord(A,N/N') \ll \left( \frac{N}{N'} \right)^{1+\epsilon}
$$ 
for all $\epsilon>0$, and we are done.
\end{proof}

%\subsection{Conclusion}
%
%Combining proposition~\ref{p:basic-bound} with
%lemma~\ref{e:count-composite} we obtain:
%\begin{prop}
%Let $N$ be square free. Then 
%\begin{equation}
%\sum_{i=1}^N |<\psi_i, \T(\V) \psi_i > |^4  
%\ll_{\A} |n|_2^{8+\epsilon} \times \frac{N 3^{\omega(N)} }{\ord(A,N)^2}
%\end{equation}
%\end{prop}

\subsection{Conclusion}

\begin{prop}
There exists a density-one sequence $S$ of integers 
such that if $n \neq 0$ and $N \in S$ then
$$
\sum_{i=1}^N 
|\langle \T(\V) \psi_i,  \psi_i \rangle |^4  \ll 
|n|_2^{8+\epsilon} \frac{N(\log N)^{14}}{\ord(A,N)^2} 
$$
\end{prop}
\begin{proof}
%Let $\gamma$ and $S$ be as in
%theorem~\ref{t:large-order-for-most-N}. 

Let $S$ be the set of integers of the form 
$N= d s^2$, where $d$ is square free, $s \leq \log
N$, and $\omega(N) \leq 3/2 \log \log N$. By
Lemmas~\ref{l:small-square-part} and \ref{l:erdos-kac}, $S$ has
density one.   

For $N=ds^2\in S$, we wish to bound the number of solutions to
\begin{equation}
\label{e:mod-N-eq}
\V ( \A^i-\A^j + \A^k - \A^l) = 0 \bmod N, \quad 
1 \leq i,j,k,l \leq \ord(A,N) 
\end{equation}
Since $N$ is not  square free we cannot apply
lemma~\ref{l:count-composite} directly. For $N=ds^2$, $d$ square-free,
we further decompose $d=d_1 \gcd(d,s)$, so that $d_1$ and 
$N/d_1 = \gcd(d,s)s^2$ are coprime.

Given $t \in \Z$ there are exactly 
$\frac{\ord(A,N)}{\ord(A,d_1)}$ solutions to $\A^i \equiv \A^t \mod d_1$
if $i \in \Z \cap [1,\ord(A,N)]$. Thus, a solution of
\begin{equation}
\label{e:mod-d-eq}
\V ( \A^i-\A^j + \A^k - \A^l) = 0 \bmod d_1, \quad 
1 \leq i,j,k,l \leq \ord(A,d_1) 
\end{equation}
lifts to at most $(\ord(A,N)/\ord(A,d_1))^4$ solutions for which $1
\leq i,j,k,l \leq \ord(A,N)$. This, together with
lemma~\ref{l:count-composite} applied to  (\ref{e:mod-d-eq}) gives
there are at most
$$
\left( \frac{\ord(A,N)}{\ord(A,d_1)} \right)^4 
|n|_2^{8+\epsilon} 
3^{\omega(d_1)} \ord(A,d_1)^2  
$$
solutions to (\ref{e:mod-N-eq}). 

Clearly  
$\omega(d_1) \leq \omega(N)$, $\ord(A,d_1)  \leq \ord(A,N)$, and since
$d_1$, $N/d_1$ are coprime, with $N/d_1\leq s^3$, we have 
$$
\frac{\ord(A,N)}{\ord(A,d_1)}  \leq  \ord(A,\frac N{d_1}) 
\ll (\frac N{d_1})^{1+\epsilon} \leq s^{3(1+\epsilon)}
$$
for all $\epsilon>0$ (by \eqref{no. Torus(N)}). 
Hence the number $\nu(N,n)$ of solutions of \eqref{e:mod-N-eq} is bounded by
\begin{equation}\label{nu(N)}
\nu(N,n)\ll |n|_2^{8+\epsilon}s^{12+\epsilon}3^{\omega(N)}\ord(A,N)^2\,.
\end{equation}

Thus we find that for $N\in S$ the number of
solutions of \eqref{e:mod-N-eq} is bounded by 
$$
|n|_2^{8+\epsilon} 
\left( \log N \right)^{12+\epsilon}
3^{3/2\log \log N} \ord(A,N)^2 \ll 
|n|_2^{8+\epsilon} (\log N)^{14}\ord(A,N)^2 
$$
and consequently we see from  Proposition~\ref{p:basic-bound} that 
\begin{equation*}
\sum_{i=1}^N 
|\langle \T(\V) \psi_i, \psi_i \rangle |^4    \leq  
|n|_2^{8+\epsilon} 
\frac{N (\log N)^{14}}{\ord(A,N)^2} 
\end{equation*}
as required. 
\end{proof}

By a routine argument (see \cite{KR}) we get:
\begin{cor}
There is a density one sequence of integers $N$ so that for all 
observables $f\in C^\infty(\TT)$, we have 
$$
\sum_{j=1}^N |\langle \OPN(f)\psi_j,\psi_j \rangle -\int_{\TT} f |^4
\ll_f \frac{N (\log N)^{14}}{\ord(A,N)^2} 
$$
\end{cor}

This reduces the proof of Theorem~\ref{main thm} to showing that for a
sequence of density one of  integers, $\ord(A,N)$ grows faster than 
$N^{1/2}(\log N)^7$ as $N\to\infty$. We will do this in
Section~\ref{ae section}  
(Theorem~\ref{t:large-order-for-most-N}).

\newpage

\section{Relating the order of $A$ modulo integers to the order modulo primes}
\label{sec:L(N)}

Our goal in this section is to show (Proposition~\ref{prop L}) 
that for a set of density one of
integers $N$, $\ord(A,N)$ is not much smaller than the product of
$\ord(A,p)$ over prime divisors $p$ of $N$. 

\subsection{}

For a set of positive integers $\mathcal M = \{m_1,\dots, m_k\}$, define 
$$
\LL(\mathcal M)  = \frac{\prod_{j=1}^k m_j}{\lcm\{m_1,\dots,m_k\}}
$$
Then $\LL(\mathcal M)$ is a positive integer, $\LL(\{m\})=1$ and
$\LL(\{m_1,m_2\}) = \gcd(m_1,m_2)$.  

 From the definition,  a prime $\ell$ divides $\LL(m_1,\dots,m_k)$ 
if and only if there are two distinct indices $i\neq j$ so that
$\ell$ divides both $m_i$ and $m_j$.

\begin{lem}\label{lemma LL}
Let $\mathcal M = \{m_1,\dots, m_k\}$, $\mathcal N = \{n_1,\dots, n_k\}$ and
suppose that $m_j \mid n_j$, $1\leq j\leq k$. Then  $\LL(\mathcal M)$ 
divides $\LL(\mathcal N)$. 
In particular, 
$$
\lcm\{m_1,\dots ,m_k\} \geq \frac{\prod_j m_j}{\LL(\mathcal N)} \;.
$$
\end{lem}
\begin{proof}
 Factor $m_j = \prod_i p_i^{\alpha_{ij}}$, 
$n_j =  \prod_i p_i^{\alpha_{ij}+\beta_{ij}}$ with 
$\alpha_{ij},\beta_{ij} \geq 0$. Then 
$\LL(\mathcal M) =\prod_i p_i^{\mu_i}$, $\LL(\mathcal N) =\prod_i p_i^{\nu_i}$
where 
\begin{equation*}
\begin{split}
\mu_i &= \sum_{j=1}^k\alpha_{ij}-\max_{1\leq j \leq k}\alpha_{ij},\\
\nu_i& = \sum_{j=1}^k(\alpha_{ij}+\beta_{ij})- 
\max_{1\leq j \leq k}(\alpha_{ij}+\beta_{ij}) \;.
\end{split}
\end{equation*}
Thus the Lemma reduces to the following easily verified inequality: For any
non-negative reals $a_j,b_j\geq 0$, $1\leq j \leq k$, we have 
$$
\sum_{j} a_j-\max_j a_j \leq \sum_j (a_j+b_j) -\max_j \{a_j+b_j\} \;.
$$
\end{proof}

\subsection{}
We need to apply these considerations to bounding $\ord(A,N)$.  
Given an integer $N$, we will write $N=ds^2$ with $d$ square-free, and
further decompose $d=d_0\gcd(d,D_A)$, so that $d_0=d_0(N)$ is
square-free and co-prime to $D_A$. 

Now define 
\begin{equation} \label{def L(N)}
L(N) = \LL(\{p-\chi(p): p\mid d_0(N) \})
\end{equation}
Since $d_0 \mid N$, we have 
$$
\ord(A,N) \geq \ord(A,d_0)  =\lcm(\{\ord(A,p):p\mid d_0\})
$$
Moreover, for $p\mid d_0$ we have $\ord(A,p) \mid p-\chi(p)$ 
and so by Lemma~\ref{lemma LL} we find
$$
\lcm(\{\ord(A,p):p\mid d_0\}) \geq \frac{\prod_{p\mid d_0}
\ord(A,P)}{L(N)}
$$
and thus 
\begin{equation}\label{bounding ord in terms of L}
\ord(A,N) \geq \frac {\prod_{p\mid d_0}\ord(A,P)}{L(N)}
\end{equation}
We will show (Proposition~\ref{prop L}) that for almost all $N\leq x$,
we have  $L(N)\leq \exp(3(\log\log x)^4)$ and consequently we get as
the main result of this section:
\begin{prop}\label{prop L}
For almost all $N\leq x$, 
$$
\ord(A,N)\geq \frac{\prod_{p\mid d_0} \ord(A,p)}
{ \exp(3(\log\log x)^4)}
$$
where $d_0$ is given by writing $N=ds^2$, with $d=d_0\gcd(d,D_A)$
square-free.  
\end{prop}

\subsection{}
 For $x\gg 1$, we set $z=z(x) = (\log\log x)^3$. 
We say that an integer is {\em $z$-smooth} if it has no prime divisors
larger than $z$.  
\begin{lem}  \label{lem:smooth}
For almost all $N\leq x$ (that is for all but $O(x/\log\log x)$), 
$L(N)$ is $z$-smooth. 
\end{lem}
\begin{proof}
Suppose that $L(N)$ is divisible by a prime $\ell>z$. From the
definition of $L(N)$, this implies that there are 
two distinct prime divisors $q_1$, $q_2$ of $d_0(N)$ so that $\ell$
divides $q_i-\chi(q_i)$, $i=1,2$. In particular, $\ell \leq x^{1/2}$. 
Thus we find two distinct primes such that 
\begin{equation}\label{condition1}
q_1q_2 \mid N \quad \text{ and } \qquad q_i = \pm 1 \mod \ell, \quad i=1,2
\end{equation}
For fixed $q_1$, $q_2$ the number of $N\leq x$ divisible by $q_1q_2$
is $[x/q_1q_2]$. Thus for fixed $\ell$, the number of $N\leq x$
satisfying \eqref{condition1} is at most 
$$
\sum_{q_1,q_2=\pm 1 \mod \ell}\frac x{q_1 q_2} \leq 
x\left(\sum_{q=\pm 1 \mod \ell} \frac 1q \right)^2
$$
By Brun-Titchmarsh (Lemma~\ref{l:merten} - recall $\ell \leq x^{1/2}$), 
this is bounded (up to constant factor) by $x(\log\log
x/\ell)^2$. Summing over all 
primes $\ell>z$, we find that the number of integers $N\leq x$ such
that $L(N)$ is divisible by some prime $\ell>z$ is at most 
$$
x(\log\log x)^2 \sum_{\ell >z} \frac 1{\ell^2} \ll 
\frac{x(\log\log x)^2} z\ll \frac x{\log\log x}
$$
\end{proof}

\begin{prop}\label{prop:L(N)}
For almost all integers $N\leq x$ we have 
$$
L(N) \leq \exp (3(\log\log x)^4) \;.
$$ 
\end{prop}
\begin{proof}
By Lemma~\ref{lem:smooth} we may assume that $L(N)$ is $z$-smooth, with
$z=(\log\log x)^3$. For $p\mid d_0(N)$,  write the $z$-smooth part of
$p-\chi(p)$ as $f_p s_p^2$, with $f_p$ square-free. 
Set 
$$
S_N=\max_{p\mid d_0} s_p\;.
$$ 
Note that since $f_p$ is square-free and $z$-smooth, it divides the
product of all primes $q\leq z$. Thus for $z \gg 1$ we have:  
$$
f_p\leq \prod_{q\leq z} q \leq e^{3 z/2} \;.
$$

Since $L(N)$ is $z$-smooth and divides $\prod_{p\mid d_0} (p-\chi(p))$,
it also divides the product $\prod_{p\mid d_0} f_p s_p^2$. 
Thus 
$$
L(N)\leq \prod_{p\mid d_0} f_ps_p^2\leq 
\prod_{p\mid d_0} e^{3z/2}S^2  \leq (e^{\frac 32z}S^2)^{\omega(N)}
%\prod_{p\mid d_0} e^{3z/2}S^2  = (e^{\frac 32z}S^2)^{\omega(N)}
$$
or 
\begin{equation}\label{L-inequality}
\frac{\log L(N)}{\omega(N)} - \frac 32 z \leq \log S_N^2
\end{equation}

Now for almost all $N\leq x$ we have (Lemma~\ref{l:erdos-kac})
\begin{equation}\label{small omega}
\omega(N) < \frac 32 \log\log x
\end{equation}
and so  by \eqref{L-inequality} if $L(N)$ is large, so is
$S_N$. Specifically, if 
$\log L(N) >3z\log\log x=3(\log\log x)^4$ then by \eqref{L-inequality},
\eqref{small omega},
we find 
$$
\log S_N^2 > z/2 = (\log\log x)^3/2
$$
We will show that this fails for almost all $N\leq x$ and thus prove
the Proposition. 

To estimate the number of $N\leq x$ for which 
$\log S_N^2 >z/2=(\log\log x)^3/2$,  
recall that by the definition of $S_N$ there is some prime $q$
dividing $d_0$ (and hence dividing $N$) so that the $z$-smooth part of
$q-\chi(q)$ is $f_q s_q^2$ and $S_N=s_q$ 
(in particular if $N\leq x$ then $S_N\leq x^{1/2}$).  
Thus there is a prime $q\mid N$ for which $q=\pm 1 \mod S^2$.

Given $q$ there are at most $[x/q]$ integers $N\leq x$ divisible by
$q$, and hence the total number of $N\leq x$ with $\log S_N^2 > z/2$ is at
most 
$$
\sum_{\exp(z/4)<S<x^{1/2}} \sum_{\substack{q=\pm 1 \mod S^2\\q\leq x}}
\frac xq
$$
By Lemma~\ref{l:merten} we have for fixed  $S<x^{1/2}$ 
$$
\sum_{\substack{q=\pm 1 \mod S^2\\q\leq x}} \frac xq
\ll \frac{x\log\log x}{S^2}
$$
and summing over $S>e^{z/4}$ gives at most 
$$
x\log\log x\sum_{S>\exp(z/4)} \frac 1{S^2} \ll 
\frac{x\log\log x}{\exp(z/4)} 
$$
Thus the number of $N\leq x$ for which 
$\log S_N^2 >z/2=(\log\log x)^3/2$ is at most 
$$
\frac {x\log\log x }{e^{z/4}} \ll x\log\log
x\exp(-\frac 14 (\log\log x)^3) = o(x)
$$
and we are done.
\end{proof}

\newpage

\section{Large order for primes}\label{Goldfeld section}

In this section we show that $\ord(A,p)$ is large for a positive
proportion of primes. Our main result here is:
\begin{thm}\label{order:thm}
Let $1/2<\eta<3/5$.  
Then the number of primes $p\leq x$ for which the order of the cat map
modulo $p$ satisfies $\ord(A,p)>x^\eta$ is at least  
$c(\eta)\pi(x) +o(\pi(x))$, where  
\begin{equation}\label{c(eta)}
c(\eta)  = \frac{ 3-5\eta}{2(1-\eta)}, \quad 1/2<\eta<3/5 \,.
\end{equation}
\end{thm} 

We first observe (following Hooley~\cite{Hooley}):  
\begin{lem}\label{small order:lem} 
The number of primes for which $\ord(A,p)\leq y$ is $\ll y^2$. 
\end{lem}
\begin{proof}
If $\ord(A,p)=k\leq y$ then $A^k=I \mod p$ and so $p\mid
\det(A^k-I)$. Thus the number of such primes is bounded by the total number
of prime divisors of the integers $\det(A^k-I)$, $k\leq y$, that is by 
$$
\sum_{k\leq y} \omega(\det(A^k-I))
$$
where $\omega(n)$ is the number of prime factors of $n$. 
Now trivially $\omega(n)\leq \log |n|$, and $|\det(A^k-I)| \sim \unit^k$
where $\unit>1$ is the largest eigenvalue of $A$. Thus we get a bound
for the number of primes as above of 
$$
\sum_{k\leq y} \omega(\det(A^k-I))\ll \sum_{k\leq y} k\ll y^2
$$
as required. 
\end{proof}

%Let $\chi=\chi_A$ be the quadratic character associated to the matrix $A$:   
%For $p$ prime $\chi(p)=1,0,-1$ according to whether $p$ splits,
%ramifies or is inert in the real quadratic field $K$ which contains the
%eigenvalues of $A$. 
%Thus $\ord(A,p)$ divides $p-\chi(p)$ for all primes $p$ not dividing 
%$D_A=4(\tr(A)^2-4)$. 

For $\eta\geq 1/2$, let $P_\eta(x)$ be the set  of primes $p\leq x$
for which there is 
a prime $q>x^\eta$, with $q\mid p-\chi(p)$. 
The main tool for proving Theorem~\ref{order:thm} is: 
\begin{prop}\label{goldfeld:prop}
For $1/2<\eta<3/5$ we have
$$
\#P_\eta(x)\geq c(\eta)\pi(x)\left(1+o(1)\right)
$$ 
with $c(\eta)>0$ given by \eqref{c(eta)}. 
\end{prop}
%\footnote{Currently it is known that $c(\eta)>0$ at
%least for  $1/2\leq \eta\leq 2/3$ \cite{}. } 

Theorem~\ref{order:thm} follows from
Proposition~\ref{goldfeld:prop} and the following observation: 
For all but $o(\pi(x))$ of the primes of $P_\eta(x)$ we have
$\ord(A,p)>x^\eta$.   
Indeed, for $p\nmid \ds$, $\ord(A,p)$ divides $p-\chi(p)$.  
For $p\in P_\eta(x)$, if $\ord(A,p)$ is not divisible by the large
factor $q>x^\eta$ of 
$p-\chi(p)$ then it divides $\frac{p-\chi(p)}q<x^{1-\eta}$ and so
$\ord(A,p)$ is smaller than $y=x^{1-\eta}$; the number of such primes
is by Lemma~\ref{small order:lem} at most $O(x^{2(1-\eta)}) = o(\pi(x))$ since
$\eta>1/2$.  Thus for all but $o(\pi(x))$ of the primes in
$P_\eta(x)$, we have  $q\mid \ord(A,p)$ and so for these primes 
$\ord(A,p)\geq q>x^\eta $.

%\newpage
%\appendix 
%\section{Goldfeld's theorem}\label{Goldfeldapp}

\subsection{Proof of Proposition~\ref{goldfeld:prop}} 
The proof of Proposition~\ref{goldfeld:prop} is a modification
of a theorem due to Goldfeld \cite{Goldfeld} from the case
of primes $p$ for which $p+a$ has a large prime factor for fixed $a$,
to the case when $a$ is allowed to vary with $p$ 
in a bounded fashion,
depending on a fixed set congruence conditions. 
%We will give the details in Section~\ref{Goldfeldapp}.  
%Let 
%$$
%P_\eta(x)=\{ p\leq x: \exists \text{ prime  }  q\mid p-\chi(p), 
%\quad q>x^\eta \}
%$$

The idea is as follows:  
By quadratic reciprocity, $\chi(p)$ only depends on the residue of $p$
modulo $\ds=4(\tr(A)^2-4)$.
Thus the number of primes in $P_\eta(x)$ is the sum over all
invertible residues $a\mod \ds$ of the number of
primes in 
$$ 
 P_\eta(x;\ds,a) = \{p\in P_\eta(x): p=a\mod \ds \}
$$
We will show 
\begin{equation}\label{goldfeld:a}
\# P_\eta(x;\ds,a) \simgeq \frac{c(\eta)}{\phi(\ds)}\pi(x) %(1+o(1))
\end{equation}
where $c(\eta)$ is given by \eqref{c(eta)}. 
Summing \eqref{goldfeld:a} over all invertible residues $a\mod \ds$ will give
Proposition~\ref{goldfeld:prop}.

%\subsubsection{Beginning} 
\subsubsection{} 
As in \cite{Goldfeld}, we consider the sum 
$$
S_a(x) = \sum_{\substack{m\leq x \\(m,\ds )=1}} 
\sum_{\substack{p\leq x \\ p=a\bmod \ds \\m\mid p-\chi(a)}}
\Lambda(m)
$$
and more generally for $y_1<y_2\leq x$, we set 
$$
S_a(y_1,y_2;x) =\sum_{\substack{y_1<m\leq y_2 \\(m,\ds )=1}} 
\sum_{\substack{p\leq x \\ p=a\bmod \ds \\m\mid p-\chi(a)}}
\Lambda(m) 
$$
This is the weighted sum over prime powers $m\in (y_1,y_2]$, coprime to
$\ds$, of the number of primes $p\leq x$, $p=a\mod m$ with 
$m\mid p-\chi(p)$.  

If $(m, \ds)=1$ then by the Chinese Remainder Theorem, there is a
unique $a_m \mod m\ds$ so that 
\begin{align*}
a_m&=\chi(a)\mod m \\
a_m&=a\mod \ds \,.
\end{align*}
Then we have  
$$
S_a(y_1,y_2;x) =\sum_{\substack{y_1<m\leq y_2 \\(m,\ds )=1}}\Lambda(m) 
\pi(x;m\ds,a_m) \,.
$$

\subsubsection{Prime powers} 
Let us first see that the contribution of proper prime powers $m=q^k$,
$k>1$, to $S_a(y_1,y_2;x)$ is at most $O(x/\log x)$, which will allow
us to ignore their contribution: 
Indeed, this contribution is bounded by 
$$
\sum_{\substack {q^k<x\\k>1}} \log q \cdot \pi(x;q^k\ds,a_{q^k}) \leq
\left(\sum_{\substack {q^k<x^{3/4} \\ k>1}} + 
\sum_{\substack {x^{3/4}\leq q^k<x \\ k>1}}\right) 
\log q \cdot \pi(x;q^k\ds,a_{q^k}) \,.
$$
By Brun-Titchmarsh \eqref{BT}, 
%if $q^k<x^{3/4}$ then $\pi(x;q^k\ds,a_{q^k})\ll x/q^k\log x$, 
if $q^k<x^{3/4}$ then $\pi(x;q^k\ds,a_{q^k})\ll x/(q^k\log x)$, 
so that the sum over $q^k<x^{3/4}$ is bounded by 
$$
\sum_{q^k<x^{3/4}} \log q \frac x{q^k\log x} \ll \frac x{\log x} 
$$
since 
$$
\sum_{q \text{ prime}}\sum_{k>1} \frac{\log q}{q^k} <\infty \,.
$$
As for the sum over $x^{3/4}<q^k<x$, we use the trivial bound 
$$
\pi(x;q^k\ds,a_{q^k})\ll \frac x{q^k\ds}< x^{1/4}
$$ 
(which comes from counting {\em integers} in an arithmetic progression) plus
the fact that the number of prime powers $q^k<x$ is 
$O(\log x/\log q)$. Since the primes 
contributing are no larger than $x^{1/2}$, we bound this sum by 
$$
\sum_{q<x^{1/2}} \log q \frac {\log x}{\log q}  x^{1/4}\ll x^{3/4}
$$
which is negligible.

\subsubsection{A reduction}
We reduce the study of $P_\eta(x;\ds,a)$ to that of
$S_a(x^\eta,x;x)$:
\begin{equation*}
\begin{split}
P_\eta(x;\ds,a) &= 
\sum_{\substack{x^\eta<q\leq x\\ q\nmid \ds \text{ prime }}}
\pi(x;\ds,a_q) \\ &\geq
\frac 1{\log x} \sum_{\substack{x^\eta<q\leq x\\ q\nmid \ds \text{ prime }}}
\log q\cdot \pi(x;\ds,a_q) \\
&= \frac 1{\log x}S_a(x^\eta,x;x) + O(\frac x{\log^2 x}) 
\end{split}
\end{equation*}
since the prime powers are negligible. (Also note that $q>x^{1/2}$ so
% ??? added "counted only once"
that each $p$ is counted exactly once in the first sum.) Thus in order
to prove 
\eqref{goldfeld:a}, we need  to show that for $\eta<3/5$,  
\begin{equation}\label{need:app}
S_a(x^\eta,x;x)\simgeq \frac{3-5\eta}{2(1-\eta)}\frac x{\phi(\ds)} \;.%+ o(x)
\end{equation}

\subsubsection{A division}
We write
$$
S_a(x)=S_a(1,\frac {x^{1/2}}{\log^c x};x) + 
S_a(\frac {x^{1/2}}{\log^c x},x^\eta;x) +  S_a(x^\eta,x;x) 
$$
with $c>1$ to be determined later. We will show 
\begin{eqnarray}
S_a(x)&\sim& \frac {x}{\phi(\ds)} \label{main m}\\
S_a(1,\frac {x^{1/2}}{\log^c x};x)&\sim& \frac 12 \frac {x}{\phi(\ds)}
\label{small m}
\\
S_a(\frac {x^{1/2}}{\log^c x},x^\eta;x) &\simleq & 
\frac{2\eta-1}{1-\eta}\frac {x}{\phi(\ds)} %\left( 1+o(1) \right) 
\label{mid m}
\end{eqnarray}
which will give \eqref{need:app} and hence our proposition. 

%\subsubsection{} 
%To evaluate $S_a(x)$, switch the order of summation and use the
%identity $\sum_{d\mid n}\Lambda(d) = \log n$ to get 
%\begin{equation*}
%\begin{split}
%S_a(x) &= \sum_{\substack{p\leq x \\ p=a\bmod \ds}}
%\sum_{\substack{m\mid p-\chi(a) \\(m,\ds )=1}}
%\Lambda(m) \\ 
%&= \sum_{\substack{p\leq x \\ p=a\bmod \ds}} \log(p-\chi(a)) 
%\sim \frac{x}{\phi(\ds)} \,.
%\end{split}
%\end{equation*}

\subsubsection{} 
To show $S_a(x)\sim x/\phi(\ds)$,  we first write  $S_a(x)$ as 
$$
\sum_{\substack{m\leq x \\(m,\ds )=1}}\Lambda(m) 
\sum_{\substack{p\leq x \\ p=a\bmod \ds \\m\mid p-\chi(a)}}1 = 
\left(\sum_{m\leq x} - \sum_{\substack{m\leq x\\(m,\ds)\neq 1}}\right)
\Lambda(m) 
\sum_{\substack{p\leq x \\ p=a\bmod \ds \\m\mid p-\chi(a)}}1
$$ 
To evaluate the sum over all $m\leq x$, we switch the order of 
summation and use the identity $\sum_{d\mid n}\Lambda(d) = \log n$ to get 
\begin{equation*}
\begin{split}
\sum_{m\leq x}\Lambda(m) \sum_{\substack{p\leq x \\ p=a\bmod \ds \\m\mid p-\chi(a)}}1
&= \sum_{\substack{p\leq x \\ p=a\bmod \ds}}
\sum_{m\mid p-\chi(a) } \Lambda(m) \\ 
&= \sum_{\substack{p\leq x \\ p=a\bmod \ds}} \log(p-\chi(a)) 
\sim \frac{x}{\phi(\ds)} \,.
\end{split}
\end{equation*}

To estimate the sum over  prime powers $m\leq x$, with $\gcd(m,\ds)\neq 1$,  
note that since the sum is only over the powers of the primes $q$ 
dividing $\ds$, it suffices to treat each such prime separately. 
We will show that each contributes at most $O_q(x/\log x)$ and thus prove 
\eqref{main m}. 

Indeed, the contribution of such a prime is 
\begin{equation*} 
\begin{split}
\log q \sum_{\substack{ k\geq 1\\q^k\leq x}}  
\sum_{\substack{p\leq x \\ p=a\bmod \ds \\q^k\mid p-\chi(a)}}1 &\leq 
\log q \sum_{q^k\leq x} \sum_{\substack{ p\leq x\\ q^k\mid p-\chi(a) }}1 \\
& \leq \log q \sum_{q^k\leq x} \pi(x;q^k,\pm 1) \;.
\end{split}
\end{equation*}
The contributing exponents  $k$ consist of those ("small" $k$'s) with 
$q^k\leq x/e$ and at most two "large" values of $k$ for which 
$x/e<q^k\leq x$.  
The  contribution of the "large" exponents  can be shown to be at most $O(1)$ 
by noting that $\pi(x;q^k,\pm 1)$ is at most the number of {\em integers} 
$n\leq x$ congruent to $\pm 1$ modulo $q^k$, which is at most 
$x/q^k +1 = O(1)$.  

For the "small" exponents ($k\geq 1$ such that  $q^k\leq x/e$), 
we use the Brun-Titchmarsh theorem  \eqref{BT} to bound 
$$
\pi(x;q^k,\pm 1) < \frac 2{1-q^{-1}} \frac {x/q^k}{\log x/q^k}
$$
and so the sum over all $k\geq 1$ with $q^k\leq x/e$  is at most 
$$
\log q \frac 2{1-q^{-1}} \sum_{q^k\leq x/e} \frac {x/q^k}{\log x/q^k} \;.
$$
In the range $q\leq q^k\leq x/e$, the function 
$k\mapsto  \frac {x/q^k}{\log x/q^k}$ is decreasing and so the sum over 
$1\leq k\leq \log(x/e)/\log q$ 
is bounded by the integral 
$$
\int_0^{\log (x/e)/\log q}\frac {x/q^k}{\log (x/q^k)} dk  = 
\frac 1{\log q} \int_e^x \frac{dt}{\log t} 
\ll \frac 1{\log q} \frac x{\log x} 
$$
Thus the total contribution of these "small" $k$'s is at most $c_q x/\log x$. 
Summing over all prime divisors $q$ of $\ds$ gives \eqref{main m}.   

\subsubsection{} 
To evaluate $S_a(1,\frac {x^{1/2}}{\log^c x};x)$, we replace
$\pi(x;m\ds,a_m)$ by $\Li(x)/\phi(m\ds)$ and use the Bombieri-Vinogradov
theorem  to bound the error by 
$$
\sum_{m<x^{1/2}/\log^c x} \log m \max_{(b,m)=1}\left|\pi(x;m\ds,b) -
\frac{\Li(x)}{\phi(m\ds)} \right|  \ll  \log x\frac{x}{\log^2 x} \ll
\frac{x}{\log x} 
$$
($c$ was chosen to give the exponent $2$ on the RHS of \eqref{BV}). 
The main term is evaluated by (note that $\phi(m\ds) =
\phi(m)\phi(\ds)$ if $m$ and $\ds$ are coprime) 
\begin{equation*}
\begin{split}
\sum_{\substack{ m<x^{1/2}/\log^c x \\(m,\ds)=1}} 
\frac{\Lambda(m)}{\phi(m\ds)}\Li(x) &=\frac{\Li(x)}{\phi(\ds)}
\sum_{\substack{ m<x^{1/2}/\log^c x\\(m,\ds)=1}} 
\frac{\Lambda(m)}{\phi(m)} \\  
&= \frac{\Li(x)}{\phi(\ds)}\left( 
\sum_{m<x^{1/2}/\log^c x}\frac{\Lambda(m)}{\phi(m)}  + O(1) \right) \\
&\sim \frac{\Li(x)}{\phi(\ds)} \log \frac{x^{1/2}}{\log^c x} \\
&\sim \frac 12\frac{x}{\phi(\ds)}
\end{split}
\end{equation*}
as required to prove \eqref{small m}.

\subsubsection{} 
Finally we estimate  $S_a(x^{1/2}/\log^c x,x^\eta;x)$, We will use the
Brun-Titchmarsh inequality \eqref{BT} which for $m<x^\eta$, $\eta<3/5$
gives 
\begin{equation}\label{BT2}
\pi(x;m\ds,a_m)< \frac 2{1-\eta} \frac x{\phi(\ds m)\log x} \;. 
%\left(1+O(\frac 1{\log x}) \right) \,.
\end{equation}
%We will neglect the O-term above as will be justified by our estimation
%using the first term. 
We now find using \eqref{BT2} that 
\begin{equation*}
\begin{split}
S_a(\frac{x^{1/2}}{\log^c x},x^\eta;x)& < \frac 2{1-\eta} \frac x{\log x} 
\sum_{\substack{ x^{1/2}/\log^c x<m\leq x^\eta\\(m,\ds)=1}}
\frac{\Lambda(m)}{\phi(m\ds)}   \\
&= \frac 1{\phi(\ds)} \frac 2{1-\eta} \frac x{\log x} 
\left( \log x^\eta-\log \frac{x^{1/2}}{\log^c x} +O(1) \right) \\
& \sim \frac{2(\eta-1/2)}{1-\eta} \frac x{\phi(\ds)} 
\end{split}
\end{equation*}
which gives the required estimate \eqref{mid m}.

\newpage
%\input{aeinteger}

%\newpage

\section{Large order for almost all integers}\label{ae section} 

In this section we will show that for a density one subsequence of the
positive integers, the order of $\A$ is large enough to
give give uniform distribution of all eigenfunctions of $\UN(A)$. 
We will show:
\begin{thm}
\label{t:large-order-for-most-N}
There exist $\delta>0$ and a density one subset $S$ of the integers 
such that for all $N \in S$ we have 
$$
\ord(A,N) \gg N^{1/2} \exp((\log N)^{\delta}) \,.
$$
\end{thm}

Fix $1/2<\eta <3/5$. We say that a prime $p$ is {\em good} if 
$p \nmid \ds$ and $\ord(A,p)\geq p^{\eta}$. 
Let $P_G$ be the set of good primes, and let $P_G(x)$
be the set of primes in $P_G$ that does not exceed $x$. As shown in
Theorem~\ref{order:thm}, there exists $\gamma=\gamma(\eta)>0$ such that  
$$
P_G(x) \simgeq \gamma \pi(x) \, .
$$

If $p\mid \ds$ or $\ord(A,p) < p^{\eta}$ we call $p$ {\em bad}, and if
$p\mid \ds$ or $\ord(A,p) <p^{1/2}/\log p$ we call $p$ 
{\em terrible}. 
As for good primes we let
$P_B$ and $P_T$ denote the set of bad, respectively terrible,
primes (note that $P_T \subset P_B$), and by $P_B(x)$ resp. $P_T(x)$
the number of primes less than $x$ in these sets. 
Since $P_B$ is the complement of $P_G$ which has lower density
$\gamma$, we have 
\begin{equation}
\label{l:prime-density}
P_B(x) \simleq (1-\gamma)\pi(x)
\end{equation}
As for the size of $P_T$, it is immediate from Lemma~\ref{small
order:lem} that  
\begin{equation}\label{l:conv-prime-sum}
P_T(x) = O(\frac{x}{\log^2 x}) \,.
\end{equation}

Given an integer $N$ we write $N = N_G N_B$ where 
$$
N_G = \prod_{ \substack{p_i^{a_i} \parallel N \\ p_i \in P_G}}
p_i^{a_i}, \quad 
N_B =\prod_{ \substack{p_i^{a_i} \parallel N \\ p_i \in P_B}}
p_i^{a_i}.  
$$
We also let $N_T \mid  N_B$ be given by $N_T =\prod_{ \substack{p_i^{a_i}
    \parallel N \\ p_i \in P_T}} p_i^{a_i}$.

Define a set of integers $\N_G$ by $n \in \N_G$ if and only if all prime
divisors of $n$ are good, and similarly for $\N_B$ and $\N_T$. As for
primes we let 
$\N_G(x)$ (respectively $\N_B(x)$ and $\N_T(x)$) be the elements of $\N_G$
(respectively  $\N_B$ and $\N_T$) not exceeding $x$.

\begin{prop}\label{p:wirsing}
The number $N_B(x)$ of integers $N\leq x$ having all their prime factors in
$P_B$ satisfies
$$
\N_B(x) \ll \frac{x}{(\log x)^\gamma} \,. 
$$
\end{prop}
\begin{proof}
Let $b_p=1$ if $p \in P_B$ and let $b_p = 0$ if $p \in P_G$, and for
composite integers $d$ put
$b_d = \prod_{p\mid d} b_p$. Then $\N_B(x) = \sum_{n \leq x}
b_n$. Since $P_B(x) \leq (1-\gamma) \pi(x)$ the sieve of Eratosthenes
gives that $\N_B(x) = o(x)$. Indeed,
$$
N_B(x)=\#\{ n \leq x \colon p   \in P_G \Rightarrow p \nmid n\} 
= x \prod_{p \in P_G(z)}(1-1/p)   + O(\exp(z)) \,. 
$$
Putting $z=\log \log x$ and noting that $\lim_{z \rightarrow
  \infty}\prod_{p \in P_G(z)}(1-1/p) = 0$ since $\sum_{p \in P_G} 1/p
= \infty$ we obtain $N_B(x) = o(x)$.

Now following  Wirsing \cite{Wirsing}, 
we consider the smoothed sum $\int_1^x \N_B(t)\frac{dt}t$. 
By partial summation we have  
\begin{equation}
\int_1^x \N_B(t)\frac{dt}t = 
\N_B(x) \log x  - \sum_{n \leq x} b_n \log n \,.
\end{equation}
Using the identity $\log n = \sum_{d\mid n} \Lambda(d)$ we obtain:
\begin{equation}
\begin{split}
\sum_{n \leq x} b_n \log n & =
\sum_{n \leq x} b_n \left( \sum_{d\mid n} \Lambda(d) \right) = 
\sum_{d \leq x} b_d \Lambda(d) \sum_{n \leq x/d } b_n \\
& = \sum_{n \leq x} b_n  \sum_{d \leq x/n } b_d \Lambda(d) \,.
\end{split}
\end{equation}
Now, 
$$
\sum_{d \leq x/n } b_d \Lambda(d) = \sum_{p \in P_B(x/n)} \log p  +O(
(\frac{x}{n})^{1/2} \log(x/n)) \ll \frac{x}{n}
$$
by Chebyshev's bound on $\pi(x)$. 
Moreover, $\N_B(t) = o(t)$ implies that $\int_1^x
\frac{\N_B(t)}{t} dt = o(x)$. Hence
\begin{equation}
 \N_B(x) \log x + o(x) \ll \sum_{n \leq x} b_n \frac{x}{n}.
\end{equation}
However, 
\begin{equation*}
\begin{split}
\sum_{n \leq x} 
\frac{b_n}{n} & \leq \prod_{p \in P_B(x)} (1+1/p+1/p^2+ \ldots) = 
\exp( \sum_{p \in P_B(x)} ( 1/p +O(1/p^2) ) \\ 
& \ll \exp \left( (1-\gamma) \log \log x \right) =
(\log x)^{1-\gamma}
\end{split}
\end{equation*}
and thus
\begin{equation}
\N_B(x) \ll \frac{x}{\log x} (\log x)^{1-\gamma}  +
o(\frac{x}{\log x}) 
\ll \frac{x}{(\log x)^\gamma} \,.
\end{equation}
\end{proof}

\begin{cor}
\label{p:large-good-part}
We have
\begin{equation}
\label{e:small-good-part}
%\#\{ N \leq x : N_G \leq z \} \ll x \frac{\log z}{ (\log \frac{x}{z})^{\gamma}}
\#\{ N \leq x : N_G \leq \exp( (\log x)^{\gamma/2}) \} 
\ll \frac{x}{  (\log x)^{\gamma/2}} 
\end{equation}
\end{cor}
\begin{proof}
We may write $\#\{ N \leq x : N_G \leq z \}$ as 
$$
\sum_{N_G \leq z} N_B(\frac{x}{N_G}),
$$
and by Proposition~\ref{p:wirsing} we may bound this sum by
$$
\sum_{N_G \leq z} \frac{x}{N_G ( \log \frac{x}{N_G} )^{\gamma} } 
\ll
\frac{x}{( \log \frac{x}{z} )^{\gamma}}
\sum_{N_G \leq z} \frac{1}{N_G} 
\ll
\frac{x}{( \log \frac{x}{z} )^{\gamma}} \log z \,.
$$
Putting $z = \exp( (\log x)^{\gamma/2})$ we obtain the desired conclusion.
\end{proof}

We will also need to estimate the number of integers $N$ with $N_T$
large: 
\begin{lem}\label{l:small-terrible}
Let $\beta(z) = \sum_{\substack{ N \in \N_T \\ N \geq z}} 1/N$. Then:  

i) The number of integers $N\leq x$ for which $N_T \geq z$ is at most 
$x \beta(z)$. 

ii) $\lim_{z \to \infty} \beta(z) =0$.
\end{lem}
\begin{proof}
i) We have 
$$
\#\{ N \leq x : N_T \geq z \} 
\leq \sum_{N_T \geq z} \frac{x}{N_T}
= x \beta(z).
$$

ii) By \eqref{l:conv-prime-sum}, 
$$
\sum_{p \in P_T} 1/p < \infty
$$
and hence
$$
\sum_{N \in \N_T} 1/N = \prod_{p \in P_T} (1 + 1/p +1/p^2 + \ldots)
< \infty.
$$
\end{proof}

\noindent{\bf Proof of Theorem~\ref{t:large-order-for-most-N}}.  
As in section~\ref{sec:L(N)}, write $N=ds^2$ where $d$ is square free, 
$d=d_0\gcd(d,D_A)$, $D_A=4(\tr(A)^2-4)$.  
%and let $L(N)$ be given by 
%\eqref{def L(N)} (see section~\ref{sec:L(N)} for more details). 
By Proposition~\ref{prop L},  for almost all $N\leq x$ we have 
\begin{equation*}
\ord(A,N) \geq \frac{ \prod_{p\mid d_0} \ord(A,p) }{\exp(3(\log\log x)^4)}\;.
\end{equation*}
Fix $1/2<\eta<3/5$. 
Write $d_0=d_G d_B$ where $d_G$ is ``good'' and $d_B$ is ``bad''. 
By definition, if $p$ is good then $\ord(A,p) >
p^{\eta}$, hence 
$$
\prod_{p\mid d_G} \ord(A,p) \geq d_G^{\eta}.
$$
Furthermore, 
$$
\prod_{p\mid \frac{d_B}{d_T}} \ord(A,p) \geq
\prod_{p\mid \frac{d_B}{d_T}} \frac{p^{1/2}}{\log p} \geq
\left( \frac{d_B}{d_T} \right)^{1/2} \frac{1}{(\log d_B)^{\omega(d_B)}}.
$$
But trivially $\ord(A,p) \geq 1$ for $p \in P_T$, hence
\begin{equation*}
\begin{split}
\prod_{p\mid d_0} \ord(A,p) & \geq
d_G^{\eta}
\left( \frac{d_B}{d_T} \right)^{1/2} \times
\frac{1}{(\log d_B)^{\omega(d_B)}} \\
& = 
\frac{ d_G^{\eta-1/2} d^{1/2} }
{ d_T^{1/2} (\log d_B)^{\omega(d_B)} } \\
& =
\frac{ d_G^{\eta-1/2} N^{1/2} }
{ (d_T s^2)^{1/2}(\log d_B)^{\omega(d_B)} } \,.
\end{split}
\end{equation*}
Now consider $N \leq x$. By the previous results we may, without
affecting the density (i.e. for all but $o(x)$), assume that the
following holds:
\begin{align}
d_T & \leq \log x && \text{(Lemma~\ref{l:small-terrible})} \\
s & \leq \log x && \text{(Lemma~\ref{l:small-square-part})} \\
\omega(d_B) & \leq  \omega(N) \leq 2 \log \log x &&
\text{(Lemma~\ref{l:erdos-kac})} \\ 
d_G & \geq \exp( (\log x)^{\gamma/2} ) &&
\text{(Corollary~\ref{p:large-good-part})} 
\end{align}
We also use $\log d_B \leq \log N\leq \log x$. Hence
$$
\prod_{p\mid d_0} \ord(A,p) \geq 
\frac{ N^{1/2} \exp \left( (\eta - 1/2) (\log x)^{\gamma/2} \right) }
{(\log x)^{3/2 + 3/2 \log \log x}} \;.
$$
%Furthermore, by Proposition~\ref{prop:L(N)} we may assume that 
%$$
%L(N) \leq \exp \left( 3 (\log \log x)^4 \right) 
%$$
Hence by Proposition~\ref{prop L}, 
\begin{equation*}
\begin{split}
\ord(A,N) \geq & \frac{ \prod_{p\mid d_0} \ord(A,p)}{\exp(3(\log\log x)^4)} 
\\
\geq &
\frac{ N^{1/2} \exp \left( (\eta - 1/2) (\log x)^{\gamma/2} \right) }
{ \exp \left( 3  (\log \log x)^4 + ({3/2 + 3/2 \log \log x}) \log \log x
  \right) } 
\\
\gg &
N^{1/2} \exp( (\log N)^{\gamma/3}  ) \,.
\end{split}
\end{equation*}
This concludes the proof of Theorem~\ref{t:large-order-for-most-N}. \qed

\newpage

\appendix
\section{Background from prime number theory}
\subsection{} 
In this Appendix, we collect some facts  which
we will need in the rest of the paper.  
The first asserts that most integers have only small square factors:
\begin{lem}
\label{l:small-square-part} 
The number of integers $N\leq x$ which have a square factor $s^2\mid
N$ with $s>\log N$ is $o(x)$. 
\end{lem}
\begin{proof}
If $N \in [x^{1/2},x]$ then $\log N \geq 1/2 \log x$, and the
number of $N \in [x^{1/2},x]$ such that $s^2 \mid N$ for some 
$s > \log N$ is bounded by 
$$
\sum_{s \geq 1/2 \log x} \frac{x}{s^2} \ll \frac{x}{\log x} \,.
$$
Hence the number of $N \leq x$ for which $s^2 \mid N$ for some 
$s>\log N$ is $\ll \frac{x}{\log x} + x^{1/2} =o(x)$. 
\end{proof}

\subsection{}
We  will need to know that  most integers have few prime factors: 
Let $\omega(N)$ be the number of prime factors of $N$. 
As a consequence of the Hardy-Ramanujan theorem   \cite{Har-Ram} (see
\cite{Hardy79}, Theorem 431), we have:  
\begin{lem}
\label{l:erdos-kac}
  The set of $N$ such that $\omega(N) \geq 3/2 \log\log N$ has zero
  density.
\end{lem}

\subsection{}
We recall two important theorems: The first is the Brun-Titchmarsh
inequality, which we will use in the following convenient form 
\cite{MV}:  
%(see \cite{HR}, page 110, thm 3.8): 
For all $1\leq k<x$, $ (a,k)=1$ 
\begin{equation}\label{BT}
\pi(x;k,a)< \frac{2x}{\phi(k)\log \frac xk} \;.
%\left( 1+\frac 8{\log \frac xk} \right)
\end{equation}
One consequence we will need is: 
\begin{lem}
\label{l:merten}
 Let $q \leq x^{1/2}$.  Then
$$
\sum_{ \substack{ p \leq x \\ p \equiv \pm 1 \mod q }} \frac{1}{p}
\ll
\frac{\log \log x}{\phi(q)} \;.
$$
\end{lem}

The second is the Bombieri-Vinogradov theorem \cite{Bombieri} in the form: 
For every $A>0$ there is some $B>0$ so that 
\begin{equation}\label{BV}
\sum_{k\leq \frac{x^{1/2}}{(\log x)^B}} \max_{(a,k)=1} \left|
\pi(x;k,a)-\frac{\Li(x)}{\phi(k)} \right| \ll \frac x{(\log x)^A} \;.
\end{equation}

\newpage

\end{document}